\documentclass[final]{article}
\usepackage{amssymb,amsmath,amsfonts,amsthm,mathrsfs,graphicx,color}
\usepackage[colorlinks=true,linkcolor=blue]{hyperref}

\usepackage{amsmath,amssymb,amsfonts,mathrsfs,color}
\usepackage{graphicx,epsf}
\usepackage[colorlinks=true,linkcolor=blue]{hyperref} % MOD
\newtheorem{theor}{{\bf Theorem}}[section]
\newtheorem{lemma}[theor]{{\bf Lemma}}           
\newtheorem{coro}[theor]{{\bf Corollary}}
\newtheorem{prop}[theor]{{\bf Proposition}}
\newtheorem{remark}[theor]{{\bf Remark}}
\numberwithin{equation}{section}

\newcommand{\refq}[1]{~(\ref{#1})}

\newcommand{\egal}{\Longleftrightarrow}
\newcommand{\epsl}{\varepsilon} 
\newcommand{\id}{{\rm id}}
\newcommand{\calA}{\mathcal A}
\newcommand{\calD}{\mathcal D}

\newcommand{\calH}{\mathcal H}
\newcommand{\calL}{\mathcal L}
\newcommand{\calO}{\mathcal O}
\newcommand{\Om}{\Omega} 
\newcommand{\ove}{\overline}

\newcommand{\N}{\mathbb N}
\newcommand{\R}{\mathbb R}

\title%[Recovering time-dependent singular coefficients]
{Recovering time-dependent singular coefficients of the wave-equation - One Dimensional Case} 

\author{O. Poisson\thanks{Aix Marseille Universit{\'e}, I2M, UMR CNRS 6632,
 France ({\tt  olivier.poisson@univ-amu.fr}).}}

\date{\today}

\begin{document}
\baselineskip 14pt
\maketitle

%\begin{abstract}
%\end{abstract}

%%%%%%%%%%%%%%%%%%%%% Section 1 %%%%%%%%%%%%%%%%%

\section{Introduction}
 Let $\Om=]0,b[\subset\R$, $b>0$, and consider the following initial boundary value problem
\begin{eqnarray}
\label{eq.u}
 \left\{ \begin{array}{rll}
  \calL_\gamma u &=&0\quad {\rm in} \quad (0,T)\times\Om=\Om_T,\\
  u|_{x=0} &=& f(t)  \quad {\rm on} \quad (0,T), \\
  u|_{x=b} &=& 0  \quad {\rm on} \quad (0,T), \\
  u\big|_{t=0} &=& u_0  \quad {\rm on} \quad \Om,\\
  \partial_t u\big|_{t=0} &=& u_1  \quad {\rm on} \quad \Om,
\end{array} \right.
\end{eqnarray}
 where $\calL_\gamma u=\partial^2_t u - \nabla_x\cdot(\gamma \nabla_x u)$,
 $\gamma = \gamma(t,x)$ has the following properties :\\
 There exist a positive constant $k \neq 1$ and a smooth function $t\mapsto a(t)\in ]0,b[$
 such that
\begin{eqnarray}
\label{def.gamma}
 \gamma(t,x) = \left\{ 
\begin{array}{ll}
 1 \quad {\rm if} \quad x <a(t),\\
 k^2 \quad {\rm if} \quad  x\in ]a(t),b[ = D(t).
\end{array} \right.
\label{S1gammacond}
\end{eqnarray}
  We make the following assumption %\HOX{\`a remplacer par $ \|\dot a(t)\|_\infty < 1$?}
$$ {\rm (H1D)}\quad \|\dot a(t)\|_\infty < \min(1,k), $$
 where $\dot a=\frac{da}{dt}$.
 The inverse problem were are concern with is to obtain some informations on $a(\cdot)$
 and $k$, by choosing carefully the data $f$ and then measuring $\partial_x u(t,x)$ at $x=0$.

\vskip.5em
 Since the velocity of waves in $\Om\setminus D(t)$ is one, it is quite natural to
 consider the following functions. We set
\begin{eqnarray}
\label{d.xit}
 \xi(t) = t-a(t),\\
\label{d.mut}
 \mu(t)=t+a(t).
\end{eqnarray}
%
% (Thanks to Assumption (H1D), these definitions hold).
 For simplicity, and if it is unambiguous, we shall write $\xi(t)=\xi$, $\mu(t)=\mu$.\\
 If needed, we extend $a(t)$ in $\R\setminus [0,T]$ by a smooth extention, and
 so we extend $D=\{\{t\}\times(a(t),b))$, $t\in [0,T]\}$, $D^C=\{\{t\}\times (0,a(t))$,
 $t\in [0,T]\}$, $\partial D=\{(t,a(t))$, $t\in [0,T]\}$ too (with the same notation)
 by replacing $[0,T]$ by $\R$ in their definition, in such a way that
$$ \delta:=\frac12 {\rm d}(\partial D,\R\times\Om)>0,\quad |\dot a|_\infty<\min(1,k). $$
 We put
$$ t_s := \inf\{t\ge s; a(t)=t-s\},\quad t^\ast(s)=2t_s-s, \quad  s\in [0,T] .$$
\begin{remark}
\label{r.ts}
 Since $|\dot a|<1$ and $a>0$, it becomes obvious that $\{t\ge s; a(t)=t-s\} = \{t_s\}$, and
 that $s\mapsto t_s$ and $s\mapsto t^\ast(\cdot)$ are smooth and increasing.
\end{remark}
 In fact, $t_0$ is the necessary time delay to have the first information on $D(t)$, and
 $t_s$ is the same, but with initial time at $t=s$.
 We set
$$ \mu_0 := t_0+a(t_0)=2t_0.$$
\begin{remark}
\label{r.ts-1}
 We obviously have $\mu(t_s) = t^*(s)$ and $\xi(t_s)=s$. Hence $\mu=t^*\circ\xi$ and
 $\xi^{-1}(\cdot)=t_{(\cdot)}$.
\end{remark}
 We also define the coefficient of reflexion/transmition by
\begin{eqnarray}
\label{d.alpha}
 \alpha(t) &:=& \frac{1-k + (k-\frac1k)\dot a(t)}{1+k+(k-\frac1k)\dot a(t)}
  = \left( \frac{1-k}{1+k} \right) \left(\frac{1- (1+\frac1k)\dot a(t)}{1+(1-\frac1k)\dot a(t)}\right) ,\\
\label{d.beta}
 \beta(t) &:=& \frac2{1+k+(k-\frac1k)\dot a(t)}.
\end{eqnarray}
 Thanks to (H1D), the functions $\alpha$ and $\beta$ are well-defined in $[0,T]$.
 We shall deal with data and measurements as functions in the usual Sobolev space $H^s(I)$,
 where $s\in\R$ and $I\subset\R$ is an non empty open interval.
 If $s\in (0,1)$ it can be defined by
$$
 H^s(I) = \Big\{q\in L^2(I);\; \int\!\!\!\int_{I\times I} \frac{|q(x)-q(y)|^2}{|x-y|^{1+2s}} dx\, dy
  <\infty \Big\}, \quad 0<s<1.
$$
 Our main result is the following

\begin{theor}
\label{t.1}
 Assume that $(u_0,u_1)\in H^{r_0}(\Om)\times H^{r_0-1}(\Om)$ for some $r_0\in (0,\frac12)$.
 Fix $f\in L^2(-\infty,T)$ such that
\begin{enumerate}
\item $f|_{(-\infty,0)} \in H^{r_0}(-\infty,0)$;
\item $f|_{(0,t)} \in H^{r_0(1-t'/T)}((0,t))$ for $0< t<t'\le T$;
\item $f|_{(0,t')} \not\in H^{r_0(1-t/T)}((0,t'))$ for $0\le t<t'\le T$.
\end{enumerate}
%
% Let $f(t)=G(t/T+r_0R)|_{[0,T]}$, where $G$ comes from Lemma \ref{l.fctG} with some $R>2T$ 
 Then, the following statements hold.\\
1) There exists a unique solution $u$ of\refq{eq.u} in $L^2(\Om_T)$.\\
2) The quantity $\partial_x u|_{x=0}$ is defined in $H^{-1}(0,T)$ by continuous extension.\\
3) The distribution $g=\partial_x u|_{x=0}+ f' \in H^{-1}(0,T)$ has the following form
$$ g = g_A+g_E,$$
 where $g_A,g_E$ satisfy  the following properties:
\begin{enumerate}
\item[(i)] $g_A(\mu) = 2\alpha(t) f'(\xi) , \quad\forall \mu \in [0,T]$.
%\begin{eqnarray*}
%\label{v.gA}
%\end{eqnarray*}
%
%
\item[(ii)] $g_A|_{(0,\mu)} \in H^{r_0(1-\tilde\xi/T)-1}(0,\mu)$ for all $\mu_0<\mu\le T$ and all $\tilde\xi>\xi$.
\item[(iii)] If $\dot a(t)\neq \frac{k}{1+k}$ then
 $g_A|_{(0,\mu)} \not\in H^{r_0(1-\tilde\xi/T)-1}(0,\mu)$, $\forall \tilde\xi<\xi$.
%
%% COMMENTs: 1) il faut T\ge\mu_0\ge\mu, car f'(\mu) intervient dans la mesure de g(\mu) .
%%                      2) \mu(T)>T car a(t)>0.
%
\item[(iv)] There exists $\epsl>0$ such that
\begin{equation}
\label{e.gE}
 g_E|_{(0,\mu)} \in H^{\epsl+r_0(1-\xi/T)-1}(0,\mu), \quad \forall \mu\in [0,T].
\end{equation}
\end{enumerate}

\end{theor}
 The main consequence of this is
\begin{coro}
\label{c.main}
 Assume that $\dot a(t)\neq \frac{k}{1+k}$ for all $t$, and
 $(u_0,u_1)\in H^{r_0}(\Om)\times H^{r_0-1}(\Om)$ for some $r_0>0$.
 Let $T>0$. We claim that:\\
1) We can know if $T\le \mu_0$ or if $T>\mu_0$.\\
2) Assume that $T>t^\ast(0)=\mu_0$. Set
$$ s^\ast := {t^\ast}^{-1}(T),\quad t_{max}:=t_{s^\ast}.$$
 Then we can recover the functions $s\mapsto t_s$, $0\le s\le s^\ast$, $t\mapsto a(t)$,
 $t_0\le t\le t_{max}$.
 The constant $k$ is the root of a second degree equation with known coefficients.
 If $\dot a\le 0$ then this equation has no more than one positive root, and so,
  we are able to reconstruct $k$.
\end{coro}
\begin{remark}
\label{r.T>mu0}
 Obviously, from Corollary \ref{c.main} and Remark \ref{r.ts}, and since $t_0=a(t_0)<b$,
 we can ensure the condition $T> \mu_0$ by choosing $T\ge 2b$.
\end{remark}
\vskip.5em
 In Theorem \ref{t.1}, the existence of such a function $f$ is ensured, thanks to the following
\begin{lemma}
\label{l.fctG}
 For all $R>0$, there exists a function $G(t)$, $0\le t\le 1$, such that
\begin{enumerate}
\item $G|_{(0,t)} \in H^{(1-t')/R}(0,t)$ for all $0<t<t'\le 1$.
\item $G|_{(0,t')} \not\in H^{(1-t)/R}(0,t')$ for all $0<t<t' \le 1$.
\end{enumerate}
\end{lemma}
\begin{remark}
 In Theorem \ref{t.1}, if $(u_0,u_1)\in H^{r_0}_0(\Om)\times H^{r_0-1}(\Om)$ for some
 $r_0\in (\frac12,1]$, and if $u_0(0)$ is known, then we can fix  $f\in L^2(0,T)$ such that
\begin{enumerate}
\item $f|_{[0,t]} \in H^{r_0(1-t/T)}([0,t])$ for $0<t\le T$;
\item $f|_{[0,t']} \not\in H^{r_0(1-t/T)}([0,t'])$ for $0<t<t'\le T$,
\end{enumerate}
 and with $f(0)=u_0(0)$. Then, the same result holds than in Theorem \ref{t.1}, but with
 $r_0\in (\frac12,1]$.\\
 If $(u_0,u_1)\in H^{r_0}(\Om)\times H^{r_0-1}(\Om)$ for some $r_0\in (\frac12,1]$,
 but if we don't know the value of $u_0(0)$, then the information is not sufficient
 (with our approach) to construct $f$ so that  the result of Theorem \ref{t.1} holds
 {\em with this value} $r_0\in (\frac12,1]$, and so, we are obliged to come back to
 the situation $(u_0,u_1)\in H^{r_1}(\Om)\times H^{r_1-1}(\Om)$, where $r_1<\frac12$.
\end{remark}
 The paper is organized as follows.  In Section \ref{s.1}, we analyse the direct problem\refq{eq.u}.
 In Section \ref{s.Ansatz} we construct an ansatz $u_A$ for\refq{eq.u} where $f$ is the function
 of Theorem \ref{t.1}.
 In Section \ref{s.3}, we first prove Corollary \ref{c.main}, then Theorem \ref{t.1}.
 In particular, we analyse the error $u_E=u-u_A$.

%%%%%%%%%%%%%%%%%%%%%%%%%%%%%%%%
%%%%%%%%%%  DIRECT PROBLEM  %%%%%%%%%%%%%%
%%%%%%%%%%%%%%%%%%%%%%%%%%%%%%%%
\section{Study of the direct problem}
\label{s.1}
\subsection{Notations}
\label{ss.1}
 We denote by $(|)$ the usual scalar product in $L^2(\Om;dx)$, by $(|)_H$ the scalar product
 in a Hilbert space $H$,  by $<;>_{H^*\times H}$ the duality product between a Hilbert space
 $H$ and its dual space $H^*$, by $<;>$ the duality product in $\calD'(\Om_T)\times \calD(\Om_T)$
 or in $\calD'(0,T)\times \calD(0,T)$.
 We put $\calH^1=L^2(0,T;H^1_0(\Om))$, $\calH^{-1}=L^2(0,T;H^{-1}(\Om))={\calH^1}^*$,
 $W=\{v\in \calH^{-1} ;\; \partial_t v\in \calH^{-1}\}$ with obvious norms.
 We denote
$$ E^r=H^r(\Om)\times H^{r-1}(\Om)\times H^r(0,T) , $$
 and
$$
 E^r_0= \left\{ \begin{array}{lc}
 \{(u_0,u_1,f)\in E^r ; \; u_0(0)=f(0), u_0(b)=0\},  & \frac12< r\le 1, \\
% \{(u_0,u_1,f)\in E^r ;\; u_0(0)"="f(0), u_0(b)"="0\},  & r=\frac12, \\
  E^r, & 0\le r < \frac12 .
\end{array}\right.
$$
 (For $r=\frac12$ we could set $E^r_0$ as in the case $r>\frac12$, but
 the relations $u_0(0)=f(0)$ and $u_0(b)=0$ should be modified).\\
 We denote $\Om_{t_1,t_2}=(t_1,t_2)\times\Om$.

\vskip.5em
 For data $v_0,v_1,F$, let $v$ satisfying in some sense:
\begin{eqnarray}
\label{eq.v}
 \left\{ \begin{array}{rll}
 \calL_\gamma v  &=& F\quad {\rm in} \quad \Om_T,\\
 v(t,x)&=& 0  \: , \quad x\in\partial\Om,\quad  t\in(0,T), \\
 v|_{t=0} &=& v_0  \quad {\rm on} \quad \Om,\\
 \partial_t v|_{t=0} &=& v_1  \quad {\rm on} \quad \Om.
\end{array} \right.
\end{eqnarray}
 We formally define the following operators:
\begin{eqnarray*}
%\label{d.op}
 u &=& \tilde P(u_0,u_1,f), \\
 \partial_x u|_{x=0}+f' &=& \tilde Z(u_0,u_1,f), \\
 (u|_{t=s},\partial_t u|_{t=s}) &=& \tilde X(s) (u_0,u_1,f), \quad 0 \le s\le T,\\
 v &=& P(v_0,v_1,F),\\
 \partial_x v|_{x=0} &=& Z(v_0,v_1,F), \\
 (v|_{t=s},\partial_t v|_{t=s}) &=& X(s) (v_0,v_1,F),\quad 0 \le s\le T,
\end{eqnarray*}
 where $u$, $v$, are respectively solutions of\refq{eq.u},\refq{eq.v}.

\subsection{Main results}
 In this section and the one above, we state that Problems\refq{eq.u},\refq{eq.v} have
 a unique solution for adequate spaces.
\begin{lemma}
\label{l.vweak}
1. The operator $P$ is a continuous linear mapping from $H^1_0(\Om)\times L^2(\Om)
 \times (L^2(\Om_T)+W)$ into $C([0,T];H^1_0(\Om))\cap  C^1([0,T];L^2(\Om))$.\\
2. The operator $X(s)$ is continuous from $H^1_0(\Om)\times L^2(\Om)\times
  (L^2(\Om_T)+W)$ into $H^1_0(\Om)\times L^2(\Om)$, for all $s\in [0,T]$. 
%\label{d.P}
\end{lemma}
\begin{lemma}
\label{l.vveryweak}
1. The operator $P$ continuously extends as a continuous operator from $L^2(\Om)\times
 H^{-1}(\Om)\times \calH^{-1}$ into $L^2(\Om_T)$.\\
2. The operator $X(s)$ continuously extends as a continuous operator from
 $L^2(\Om)\times H^{-1}(\Om)\times \calH^{-1}$ into $L^2(\Om)\times H^{-1}(\Om)$,
 for all $s\in [0,T]$.
\end{lemma}

\begin{lemma}
\label{l.uvweak}
1. The operator $\tilde P$ is a continuous linear mapping from $E^1_0$ into
 $C([0,T]$; $H^1(\Om))\cap C^1([0,T];L^2(\Om))$, and continuously extends as a continuous
 operator from $E^0$ into $L^2(\Om_T)$.\\
2. The operator $\tilde X(s)$ is continuous from $E^1_0$ into $H^1(\Om)\times L^2(\Om)$,
 and continuously extends as a continuous operator from $E^0$ into
 $L^2(\Om)\times H^{-1}(\Om)$, for all $s\in [0,T]$.
\end{lemma}

\begin{lemma}
\label{l.Z}
 The operator $Z$ (respect., $\tilde Z$) is continuous from $H^1_0(\Om)\times L^2(\Om)\times
 L^2(\Om_T)$ (respect., $E^1_0$) into $L^2(0,T)$ and continuously extends as a continuous
 operator from $L^2(\Om)\times H^{-1}(\Om)\times \calH^{-1}$ (respect., $E^0$) into $H^{-1}(0,T)$.
\end{lemma}
\begin{lemma}
\label{l.UC}
 Let $t_1\in [0,T]$. Assume that $F\in \calH^{-1}$ has a compact support in $\calO(t_1)$. 
 Let $v=P(v_0,v_1,F)$. Then there exists a neighborhood $\tilde K$ of $K(t_1)$
 in $\ove{D^C}$ such that  $v|_{\tilde K}$ does not depend on $F$, that, is,
 if $v_0=v_1=0$, then $v|_{\tilde K}$ vanishes, and, in particular,
 supp $\partial_x v|_{x=0}\subset (\mu(t_1),T]$.
\end{lemma}

\subsection{Proofs}
% We define, for all $t\in [0,T]$, the operator $A(t) = -\nabla(\gamma(t,\cdot)\nabla\cdot)$
% with domain
%%
%$$D(A(t)) = \{q\in H^1_0(\Om); A(t)q\in L^2(\Om)\}. $$
%%
 Let us consider the familly of bilinear forms $b(t)$, $t\in\R$, defined by
$$ b(t;u,v) = \int_\Om \gamma(t,x)\nabla_x u(x)\; \nabla_x v(x) \; dx,\quad \forall u,v\in H^1(\Om).$$
 Lemma \ref{l.vweak} is a corollary of the following general theorem (proof in appendix),
 which is an extension of \cite[XV section 4] {DAU.ANA} where $\gamma$
 did not depend on the variable $t$.
\begin{theor}
\label{t.L1}
 Let $T>0$ and $\Om\subset\R^n$, $n\ge 1$, such that $H^1_0(\Om)$ is compact in $L^2(\Om)$.
 Let $\gamma(t,x)>0$ be such that $\gamma,\gamma^{-1} \in C^0([0,T];L^\infty(\Om))$,
 $\partial_t\gamma \in L^\infty(\Om_T)$.
 Let $F\in W\cup L^2(\Om_T)$ and $v_0\in H^1_0(\Om)$, $v_1\in L^2(\Om)$.
 Then, there exists a unique weak solution $v$ to\refq{eq.v}, that is,
 $v\in C([0,T];H^1_0(\Om))$, $\partial_t v\in C([0,T];L^2(\Om))$, $v|_{t=0}=v_0$,
 $\partial_t v|_{t=0}=v_1$, and
$$ \frac{d}{dt} (\partial_t v|\phi)+b(t;v(t,\cdot),\phi) = \; <F(t,\cdot);\phi>,$$
 in the sense of $\calD'(]0,T[)$, for all $\phi\in H^1_0(\Om)$.
 Moreover there exists a  constant $C$ such that
\begin{equation}
\label{in.en1}
   \|\partial_t v(t,\cdot)\|_{L^2(\Om)} + \|\partial_x v(t,\cdot)\|_{L^2(\Om)}\le C \left(\|F\|_{L^2(\Om_t)+W}
 + \|v_0\|_{H^1_0(\Om)} +  \|v_1\|_{L^2(\Om)}\right),\quad \forall t\in [0,T].
\end{equation}
\end{theor}
 Let us show that Lemma \ref{l.vveryweak} is a straightforward consequence of
 Lemma \ref{l.vweak} with the operator $P$ replaced by its adjoint $P^*$.
 Let $(v_0,v_1,F)\in L^2(\Om)\times H^{-1}(\Om)\times \calH^{-1}$.
 By the principle of duality, we can write\refq{eq.v} as
$$
 (v|g)_{L^2(\Om_T)}  \quad = \quad <v_1,w(0)>_{H^{-1}\times H^1_0}-(v_0|\partial_t w(0))
 + <F,w>_{\calH^{-1}\times \calH^1}, 
$$
 for all $g\in L^2(\Om_T)$, where we put $w=P^*(0,0,g)$.
 Consequently (thanks to Lax-Milgram theorem), Equation\refq{eq.v} admits
 a unique solution $v\in L^2(\Om_T)$, and this shows Point 1 of Lemma \ref{l.vveryweak}.
 Once again, we have
\begin{eqnarray*}
  <\partial_t v|_{t=T},f_0>_{H^{-1}\times H^1_0}-(v|_{t=T}|f_1) &=&
 <v_1,w(0)>_{H^{-1}\times H^1_0}-(v_0|\partial_t w(0)) \\
 && + <F,w>_{\calH^{-1}\times \calH^1},
\end{eqnarray*}
 for all $(f_0,f_1)\in H^1_0(\Om)\times L^2(\Om)$, where we put $w=P^*(f_0,f_1,0)$.
 This shows that $(v|_{t=T},\partial_t v|_{t=T})\in L^2(\Om)\times H^{-1}(\Om)$.
 This proves Point 2 of Lemma \ref{l.vveryweak} in the non-restrictive case $s=T$.

 Let us prove Lemma \ref{l.uvweak}. Let $\Phi(x)\in C^\infty(\R)$ with $\Phi(0)=1$
 and with support in $[0,a_m]$, where $a_m\le a(t)$ for all $t$.
 Let us consider $f\in H^1_{loc}(\R)$ first.  Set
\begin{equation}
\label{d.uin}
 u_{in}(t,x)=f(t-x)\Phi(x).
\end{equation}
 Problem\refq{eq.u} with unknown $u$ is (at least formally) equivalent to the following one:
 find $v=u-u_{in}$ satisfying\refq{eq.v} with
\begin{eqnarray}
\label{v.v01}
 \quad v_0(x) &=&u_0(x) -f(-x)\Phi(x),\quad v_1(x)=u_1(x) - f'(-x)\Phi(x), \\
\label{v.Lui}
 F(t,x) &=& -\calL_\gamma u_{in} \, (t,x)= -\calL_1 u_{in} \, (t,x)= -2f'(t-x)\Phi(x)
 + f(t-x)\Phi''(x).\hspace{-1cm}
\end{eqnarray}
 Relation\refq{v.Lui} shows that $F\in L^2(\Om_T)$. In fact, we have $F\in W$ also,
 since
$$ \partial_t F(t,x) = -2f''(t-x)\Phi(x) + f'(t-x)\Phi''(x),$$
 and, for all $\varphi\in \calD(\Om_T)$,
\begin{eqnarray*}
 <f''(t-x)\Phi(x),\varphi(t,x)> &=& <f''(t-x,)\Phi(x)\varphi(t,x)> \\
  &=&  <f'(t-x), \partial_x(\Phi(x)\varphi(t,x))> \le  C \|\varphi\|_{\calH^1}, 
\end{eqnarray*}
 which shows that $\partial_t F(t,x) \in \calH^{-1}$. Similarly, we have
$$ <f'(t-x)\Phi(x),\varphi(t,x)> = <f(t-x),\partial_x(\Phi(x)\varphi(t,x))\le C \|\varphi\|_{L^2(\Om_T)},$$
 which shows that $F \in \calH^{-1}$ if $f\in L^2_{loc}(\R)$ only. We set
$$
 R: \begin{array}{ccc} H^1_{loc}(\R) & \to & L^2(\Om_T)\cap W \\
  f & \mapsto & F \mbox{ defined by \refq{v.Lui} },
 \end{array}
$$
$$
 S: \begin{array}{ccc} E^1_0 & \to & H^1_0(\Om)\times L^2(\Om) \\
  (u_0,u_1,f) & \mapsto & (v_0,v_1) \mbox{ defined by \refq{v.v01}.}
 \end{array}
$$
 The above analysis shows that $R$ continuously extends as a continuous operator
 from $L^2_{loc}(\R)$ into $\calH^{-1}$.
 Similarly, $S$ continuously extends as a continuous operator from $E^0$ into
 $L^2(\Om) \times H^{-1}(\Om)$.
 Consequently, and since a solution to\refq{eq.u} can be written $u=v+u_{in}$ with
 $v=P(S(u_0,u_1,f),R(f))$, Point 1 of Lemma \ref{l.uvweak} is proved.
 Similarly, we prove Point 2 of Lemma \ref{l.uvweak}, since we have
 $\tilde X(s)(u_0,u_1,f)=X(s)(S(u_0,u_1,f),Rf) + (u_{in}|_{t=s},\partial_t u_{in}|_{t=s})$.

 Let us prove Lemma \ref{l.Z}.
 Let $(v_0,v_1,F)\in H^1_0(\Om)\times L^2(\Om)\times L^2(\Om_T)$.\\
 As above, for all $\varphi\in \calD(\R)$ such that $\varphi(T)=0$, there exists a unique solution
 $q=q_\varphi \in L^2(\Om_T)$ to
\begin{eqnarray}
\label{eq.dualq}
 \left\{ \begin{array}{rll}
 \calL_\gamma q  &=& 0\quad {\rm in} \quad \Om_T,\\
 (q(t,0),q(t,b)) &=& (\varphi,0)  \quad {\rm on} \quad (0,T), \\
 q|_{t=T} &=& 0  \quad {\rm on} \quad \Om,\\
 \partial_t q|_{t=T} &=& 0  \quad {\rm on} \quad \Om,
\end{array} \right.
\end{eqnarray}
 since it is a particular case of Lemma \ref{l.uvweak} with reversal time.\\
 Moreover, we have $q_\varphi\in C([0,T];H^1_0(\Om))$, $\partial_t q_\varphi\in C([0,T];L^2(\Om))$
 with
\begin{eqnarray}
\label{m0.qphi}
 \|q_\varphi |_{t=0}\|_{L^2(\Om)}+\| \partial_t q_\varphi |_{t=0}\|_{H^{-1}(\Om)}
 + \|q_\varphi\|_{L^2(\Om_T)} \le C\|\varphi\|_{L^2(0,T)}, \\
\label{m1.qphi}
 \|q_\varphi |_{t=0}\|_{H^1(\Om)}+\| \partial_t q_\varphi |_{t=0}\|_{L^2(\Om)}
 +\|q_\varphi\|_{\calH^1} + \|\partial_t q_\varphi\|_{L^2(\Om_T)} \le C\|\varphi\|_{H^1(0,T)}.
  \hspace{-1.5cm} 
\end{eqnarray}
 By the duality principle, and thanks to\refq{m0.qphi}, we have in the sense of $\calD'([0,T))$,
\begin{eqnarray}
\label{e.dxv0}
 <\partial_x v|_{x=0},\varphi> &=& - <v_0 ,\partial_t q_\varphi|_{t=0}> +
  <  v_1 ,q_\varphi|_{t=0}>  + <F,q_\varphi> \\
\label{e.dxv}
 &\le & C \left( \|v_0 |_{t=0}\|_{H^1(\Om)}+\| v_1 \|_{L^2(\Om)}+\| F \|_{L^2(\Om_T)} \right)
  \|\varphi\|_{L^2(0,T)},  \hspace{-1cm}
\end{eqnarray}
 which shows that $\partial_x v|_{x=0} \in L^2(0,T)$ and that $Z$ is a continuous mapping
 from $H^1_0(\Om)\times L^2(\Om)\times L^2(\Om_T)$ into $L^2(0,T)$.

 Now, let $(v_0,v_1,F)\in L^2(\Om)\times H^{-1}(\Om)\times \calH^{-1}$.
 Then, Relation\refq{e.dxv0} and Estimate\refq{m1.qphi} imply
\begin{eqnarray}
\label{e.dxv2}
 <\partial_x v|_{x=0},\varphi> &\le & C \left( \|v_0 |_{t=0}\|_{L^2(\Om)}+\| v_1 \|_{H^{-1}(\Om)}+
  \| F \|_{\calH^{-1}} \right) \|\varphi\|_{H^1(0,T)},  \hspace{-1cm}
\end{eqnarray}
 which shows that $\partial_x v|_{x=0} \in (H^1_T)'\subset H^{-1}(0,T)$, the dual space of
 $H^1_T=\{f\in H^1(0,T);f(T)=0\}$, and that $Z$ continuously extends as a continuous operator
 from $L^2(\Om)\times H^{-1}(\Om) \times  \calH^{-1}$ into $H^{-1}(0,T)$.\\
%
% Thus $\partial_x v|_{x=0} \in  L^2(0,T)$ with $\|\partial_x v|_{x=0}\|_{L^2(0,T)}\le C \|F\|_W$.\\
  This ends the proof of the property of $Z$ in Lemma \ref{l.Z}.
 Since $\partial_x u_{in}|_{x=0} = -f'$, we have $\tilde Z(u_0,u_1,f) = Z(S(u_0,u_1,f),Rf)$,
 and Point 2 of Lemma \ref{l.Z} is proved.

\qed

\vskip1em
 By the well-known Sobolev interpolation theory, we have also proved:
\begin{prop}
\label{p.Zinter}
The operator $P$ (respect., $\tilde P$) continuously maps $H^s(\Om)\times H^{s-1}(\Om)\times
 L^2(0,T;H^{s-1}(\Om))$ (respect., $E^s_0$) into $L^2(0,T;H^s(\Om))$, $s\in [0,1] (\setminus \frac12)$.\\
 The operator $Z$ (respect., $\tilde Z$) continuously maps $H^s(\Om)\times H^{s-1}(\Om)\times
 L^2(0,T;H^{s-1}(\Om))$ (respect., $E^s_0$) into $H^{s-1}(0,T)$, $s\in [0,1] (\setminus \frac12)$.
\end{prop}
 Proof of Lemma \ref{l.UC}. Denote $K= K(t_1)$. Notice that
 $K\cap \ove{D^C}=\{(t_1,a(t_1))\}$.
 We assume that $v_0=v_1=0$. Since supp$F\cap \ove{\Om_{t_1}}=\emptyset$,
 then, thanks to Lemma \ref{l.vveryweak} with $T$ replaced by $t_1$,
 $v$ vanishes in $\Om_{t_1}$.
 Let $K'={\rm int}\; K$ the interior of $K$. The function $v|_K\in L^2(K')$ satisfies
 the following equations:
\begin{eqnarray*}
 \partial^2_t v-\Delta_x v = 0 \quad {\rm in}\quad K' ,\\
 v(t,0) =  0 , \quad t_1<t<\mu(t_1),\\
 v|_{t=t_1}= \partial_t v|_{t=t_1}=0 \quad {\rm in}\quad (0,a(t_1)).
\end{eqnarray*}
 It is well-known that this implies $v|_{K'}=0$, and so, supp $\partial_x v|_{x=0}
 \subset [\mu(t_1),T]$. But since the support of $F$ does not touch $\partial K$,
 we similarly have $v|_{K_\epsl(t_1)}=0$, supp $\partial_x v|_{x=0} \subset [\mu(t_1)+\delta,T]$,
 for some $\epsl>0$ sufficiently small.\\
 However, let us give a more straightforward and simple proof to the fact that
 supp $\partial_x v|_{x=0} \subset [\mu(t_1)+\delta,T]$.
 Fix $\delta,\epsl>0$ such that $\mu(t_1)+\delta>\mu(t_1+\epsl)$ and
 supp $F\cap K_\epsl(t_1)=\emptyset$.
 Let $t_2\in [t_1,t_1+\epsl]$, $\varphi\in H^1_0(0,\mu(t_2))$ and set
 $w(t,x)=\varphi(t+x)$ for $t_2\le t\le \mu(t_2)$.
 Observe that $w=q_\varphi$ of\refq{eq.dualq}, but with $(0,T)$ replaced by $(t_2,\mu(t_2))$.
 In fact, supp $w\subset K(t_2)$, and so $w$ vanishes in $D\cap \Om_{t_2,\mu(t_2)}$.
 We then have, similarly to\refq{e.dxv0},
$$
 < \partial_x v|_{x=0,t_2<t<\mu(t_2)},\varphi> =
   - <v|_{t_2} ,\partial_t w|_{t_2}> +  < \partial_t v|_{t_2} , w|_{t_2}>  + <F,w> 
  =0
$$
 since  $v|_{t_2}= \partial_t v|_{t_2}=0$ and supp $F\cap$ supp $w=\emptyset$.
 Since $\varphi$ is arbitrarily chosen, this shows that  supp $\partial_x v|_{x=0}
 \cap (t_2,\mu(t_2))=\emptyset$, for all $t_2\in [t_1,t_1+\epsl]$.
 Hence,  supp $\partial_x v|_{x=0} \subset [\mu(t_1+\epsl),T]$.
 \qed.

\section{Ansatz}
\label{s.Ansatz}
\subsection{Notations}
 For $t\in [0,T]$ we put
$$
 K(t)=\{(s,x)\in [t,\mu(t)]\times \ove\Om;\; s+x\le \mu(t)\},\quad
 \calO(t)= \{(s,x)\in \ove{\Om_{t,T}}; \: s+x>\mu(t)\}.
$$
 (Notice that $K(t)\subset\ove{D^C}$ and $K(t)\cap\ove{D}=\{(t,a(t)\}$).\\
 For $\epsl>0$, $t\in [0,T]$, we put $K_\epsl(t)=\cup_{t\le s\le t+\epsl} K(s)$.

 If $q(x)$ is sufficiently smooth in $\Om$, then $[q]_t := q(a(t)+0)-q(a(t)-0)$.\\
 We write $g_1 \stackrel{s}{\simeq}g_2$ if $g_1$ or $g_2\in H^s(0,T)$ and
 $g_1-g_2\in H^{s+\epsl}(0,T)$ for some $\epsl>0$.\\
 We put $C^j_+=\{f\in C^j(\R);\; f|_{\R^-}=0\}$, $j\in\N$, which is dense in
 $L^2(\R^+)\approx\{f\in L^2(\R);\; f|_{(-\infty,0)}=0\}$.
 We consider for all $t\in [0,T]$ the formal operator $\calA(t)=-\nabla_x(\gamma(t,\cdot)\nabla_x)$
 defined from $H^1(\Om)$ into $H^{-1}(\Om)$ by duality:
$$
 <\calA(t)u,w>_{H^{-1}(\Om)\times H^1_0(\Om)}= (\gamma(t)\nabla_x u|\nabla_x w),
 \quad \forall u,w\in H^1(\Om)\times H^1_0(\Om) .
$$
 \vskip1em
 Let $f$ be a measurable function, we define the ansatz $u_A=U_A(f)$ for\refq{eq.u} as follows.
 Recall that $\xi(t)$ and $\mu(t)$ are defined by\refq{d.xit},\refq{d.mut}, and we have
\begin{eqnarray}
\label{d.xi0}
  \xi_0=t_0-a(t_0)=0,\\
\label{d.mu0}
 \mu_0=t_0+a(t_0)=2t_0.
\end{eqnarray}
 In addition, we put, for $t\in [0,T]$,
\begin{equation}
\label{d.nut}
 \nu=t-\frac{a(t)}k, \quad \nu_0= t_0-\frac{a(t_0)}k.
\end{equation}
 Thanks to Assumption (H1D), $t\mapsto \nu(t)$ is invertible.
 Recall also that the coefficient of reflexion/transmition, $\alpha$ and $\beta$, are defined
 by\refq{d.alpha},\refq{d.beta}. Note that we have
\begin{eqnarray}
\label{p0.ab}
 \alpha(t) \frac{d\mu}{d\xi} - \beta(t) \frac{d\nu}{d\xi} = -1,\\
\label{p1.ab}
 \alpha(t) + k\beta(t) = 1.
\end{eqnarray}
 We also define:
\begin{eqnarray}
\label{v.f2}
 f_2(\mu)  &=& \alpha(t) \frac{d\mu}{d\xi} f(\xi) ,\\
\label{v.f3}
 f_3(\nu) &=&  \beta(t) \frac{d\nu}{d\xi} f(\xi) .
\end{eqnarray}
 We put
\begin{eqnarray*}
\label{v.uA}
 u_A(t,x) &=&  \left\{ \begin{array}{cl}
  f(t-x) + f_2(t+x) - f_2(t-x) \Phi_\epsl(x), &  0\le t\le T ,\quad 0< x< a(t) ,\\
  \Phi_\epsl(x-b+2\epsl) f_3(t-\frac{x}k) , & 0\le t\le T \quad a(t) < x < b,
\end{array} \right. % \hspace{-2cm}
\end{eqnarray*}
 where we fix $\Phi_\epsl\in C^\infty(\R)$ so that $\Phi_\epsl(r)=1$ if $r<\frac12\epsl$,
 $\Phi_1(r)=0$ if $r>\epsl$, $0<\epsl \le \frac12 {\rm d}(\partial D,\partial\Om_T)$.
 It is clear that the linear operator $U_A:\: f\mapsto u_A$ is bounded from
 $L^2(\R)$ into $L^2(\Om_{\mu_T})$.
 \subsection{Properties of the Ansatz}
\begin{lemma}
\label{l.22}
 Let $f\in C^2(\R)$. Then we have\\
1) $u_A\in C^2([0,T];H^1(\Om))$, $u_A|_D\in C^2(\ove{D})$,
 $u_A|_{D^C}\in C^2(\ove{D^C})$.\\
2) There exists a smooth function $\tau(t)$ with support in $[t_0,\mu_0]$ such that
$$ {[}\gamma \partial_x u_A(t)]_t = \tau(t) f(\xi(t)) .$$
3) a) $u_A$ vanishes near $x=b$.\\
b)
 Let $g_A=\partial_x u_A|_{x=0}+f'$. Then $g_A(\mu)=2\alpha(t) f'(\xi)$ for $0\le \mu\le T$,
 where $t$, $\xi$, $\mu$ are related by\refq{d.xit},\refq{d.mut},\refq{d.nut}.\\
4) Put $F_A=\calL_\gamma u_A$ in the sense that
 $F_A(t,\cdot)=\frac{d^2}{dt^2}u_A(t)+\calA(t)u_A(t) \in H^{-1}(\Om)$ for all $t$,
 and $F_A \in C([0,T];H^{-1}(\Om))$.
 Then, $F_A$ can be written\\
 $F_A(t,x)=F_1(t,x) - \tau(t) f(\xi(t)) \delta_{a(t)}(x)$, where $\tau$ is smooth, and
 $F_1\in C([0,T]; L^2(\Om))$ is defined for  $0 \le t\le T$ by
\begin{eqnarray}
\label{v.F1}
 F_1(t,x) &=&  \left\{ \begin{array}{cl}
 \Phi_2(x)f_2(t-x)+\Phi_3(x)f'_2(t-x)  & 0< x< a(t),\\
 \Phi_4(x) f_3(t-\frac{x}k) + \Phi_5(x) f'_3(t-\frac{x}k)  ,  & a(t) < x < b,
\end{array} \right.  \hspace{-2em}
\end{eqnarray}
 where the functions $\Phi_j$ are smooth and independant of $f$,
 with compact support in $[\epsl/2,\epsl]$ for $j=2,3$, and in $[b-\epsl,b-\epsl/2]$
 for $j=4,5$.
\end{lemma}
 Proof. Point 1. is obvious, since we have, thanks to\refq{p0.ab},
\begin{eqnarray*}
 [u_A(t,\cdot)]_t &=& f_3(\nu(t)) -f(\xi(t)) - f_2(\mu(t))  \\
 &=& \left(\beta(t) \frac{d\nu}{d\xi} -1-   \alpha(t) \frac{d\mu}{d\xi} \right) f(\xi) =0.
\end{eqnarray*}
 Let us consider Point 2. For $0\le t\le T$ we have
\begin{eqnarray*}
\gamma \partial_x u_A(t,a(t)-0) &=& -f'(\xi)+ f'_2(\mu) = (-1+\alpha)f'(\xi) 
 + \frac{d(\alpha\frac{d\mu}{d\xi})}{d\mu} f(\xi) , \\
\gamma \partial_x u_A(t,a(t)+0) &=& -k  f'_3(\nu) = -k \beta(t) f'(\xi)
 - k\frac{d(\beta(t) \frac{d\nu}{d\xi})}
 {d\nu} f(\xi) .
\end{eqnarray*}
 Thanks to\refq{p1.ab} we get
$$ {[}\gamma \partial_x u_A(t)]_t =- \tau(t) f(\xi) ,$$
 with
$$
 \tau(t) =  - k\frac{d(\beta(t) \frac{d\nu}{d\xi})}{d\nu} -
  \frac{d(\alpha\frac{d\mu}{d\xi})}{d\mu}.
$$
 This ends Point 2.\\
 Let us consider Point 3 b), since 3 a) is obvious. For $0\le \mu\le T$ we have
\begin{eqnarray*}
 \partial_x u_A(\mu,0) &=& -f'(\mu)+2f'_2(\mu) = -f'(\mu)+2\alpha(t)f'(\xi).
\end{eqnarray*}
 This ends Point 3.\\
 Let us prove  Point 4.
 A short computation yields\refq{v.F1}.
 Thanks to Point 2, we obtain $F_A=F_1+\tau(t)f(\xi)$ in the required sense.
 This ends the proof of the lemma.
\qed

\vskip.2em
 We define the bounded operators $U_A$ : $C^2(\R)\ni f\mapsto u_A \in
 C^2([0,T];H^1(\Om))$, $T_0$: $C^2(\R)\ni f\mapsto T_0f\in C([0,T];H^{-1}(\Om))$
 such that $T_0f (t)=\tau(t)f(\xi)\delta_{a(t)}(x)$, and
 $T_1$: $C^2(\R)\ni f\mapsto T_1f =F_1\in C([0,T]; L^2(\Om))$,
 $T_A$: $C^2(\R)\ni f\mapsto T_Af = F_A\in C([0,T];H^{-1}(\Om))$.
 Notice that $T_0f(t)\in H^{-s}(\Om)$ for all $s>\frac12$, $t\in [0,T]$.
 Obviously we have the following propositions and Lemma.
\begin{prop}
\label{p.UA}
 The operator $U_A$ continuously extends as a bounded operator from
 $L^2(0,T)$ into $C([0,T];H^{-1}(\Om))$.
\end{prop}
\begin{prop}
\label{p.T0}
 The operator $T_0$ continuously extends as a bounded operator from $L^2(0,T)$
 into $L^2(0,T;H^{-s}(\Om))$, $\forall s>\frac12$.
\end{prop}
\begin{lemma}
\label{l.24}
1) The operator $T_A$ is continuous from $C^2(\R)$ into $L^2(0,T;H^{-1}(\Om))$ and,
 for all $s\in [0,\frac12)$, it extends as a continuous operator from $H^s(0,T)$ into
 $L^2(0,T;H^{s-1}(\Om))$.\\
2) The operator $G_A:\:f\mapsto \partial_x U_A(f)|_{x=0}+f'$ is continuous
 from $C^2(\R)$ into $C^0([0,T])$, and, for all $s\in [0,\frac12)$, it extends as
 a continuous operator from $H^s(0,T)$ into $H^{s-1}(0,T)$.\\
3) Let $f$ such as in Theorem \ref{t.1}, then $g_A:=G_A f$ satisfies (ii) and (iii)
 of Theorem \ref{t.1}.
\end{lemma}
 Proof of Lemma \ref{l.24}.
 Point 1). Thanks to Lemma \ref{p.T0}, it is sufficient to prove this with $T_A$
 replaced by $T_1$. Thanks to the interpolation theory, it is sufficient to prove that $T_1$
 is a bounded operator from $L^2(0,T)$ into $L^2(0,T;H^{-1}(\Om))$ and from
 $H^1_0(0,T)$ into $L^2(\Om_T)$, that is obvious. Hence Point 1) holds.
 Point 2) is obvious for the same reason.
 Point 3) is obvious, since $\alpha(t)\neq 0$ for all $t$.
\qed

\subsection{Modification of $F_1$}
 The regularity of $F_1$ is not sufficient for us, we replace it by the following one,
 $F_{\epsl,\tilde\mu}$, which is equivalent to $F_1$ in the sense of Lemma \ref{l.UC}.\\
 Let  $\tilde\mu\in[0,T]$, put $\tilde t=\mu^{-1}(\tilde\mu)$, $\tilde\xi=\xi(\tilde t)$,  
 $\tilde\nu=\nu(\tilde t)$, and consider a smooth function $\phi(\cdot;\epsl,\tilde\mu)$
 defined in $\R^2$ such that $\phi(t,x;\epsl,\tilde\mu)=1$ for $(t,x)\in \ove{\Om_{\tilde t}}
  \cup K_{\epsl/2}(\tilde t)$,  $\phi(t,x;\epsl,\tilde\mu)=0$ for $t\ge \tilde t+\epsl$ and
 $(t,x)\not\in K_\epsl(\tilde t)$.
 For $s\in[0,\frac12)$, $f\in H^s(\R)$ and $F_1=T_1(f)$ we put
\begin{eqnarray*}
 F_{\epsl,\tilde\mu}(t,x) &=&  F_1(t,x)\phi(t,x;\epsl,\tilde\mu).
\end{eqnarray*}
%
% and then $F_{A,\epsl}(t,x;\tilde\mu) = T_0(f)(t,x)+ F_{1,\epsl}(t,x;\tilde\mu)$.
 We have the two following properties.
\begin{lemma}
\label{l.suppFF1}
 For $\epsl<\delta$, the support of $F_1-F_{\epsl,\tilde\mu}$ is contained in $\calO(\tilde t)$.
\end{lemma}
 Proof. Since $F_1-F_{\epsl,\tilde\mu}=(1-\phi(\cdot;\epsl,\tilde\mu) F_1$, the support
 of $F_1-F_{\epsl,\tilde\mu}$ is contained in supp $(1-\phi(\cdot;\epsl,\tilde\mu))\cap$
 supp $F_1$. But supp $(1-\phi(\cdot;\epsl,\tilde\mu))\subset \ove{\Om_T}\setminus
 (\Om_{\tilde t}\cup {\rm int}(K_{\epsl/2}(\tilde t)) )$.
 Then the proof is done if we show that $(\tilde t,a(\tilde t))\not\in$ supp $(F_1-F_{\epsl,\tilde\mu})$.
 But, thanks to\refq{v.F1}, the support of $F_1$ is localized in $\{x\le \epsl\}\cup\{x\ge b-\epsl\}$
 that does not touch $\partial D$. 
\qed

\begin{lemma}
\label{l.pF1}
 Let $f$ as in Theorem \ref{t.1}. There exists $c>0$ and $\epsl_0>0$, independent of $f$,
 such that, for all $\epsl\in (0,\epsl_0)$, $\tilde\mu\in [0,T]$, $F_{\epsl,\tilde\mu} \in
  C([0,T];H^{r_0(1-\tilde\xi/T)+c\epsl-1}(\Om))$.
\end{lemma}
 To prove it, we use the following well-known property.
\begin{prop}
\label{p.fF}
 Let $g\in H^s(\R)$ for some $s\in [-1,0]$. Let $r\in \R^*$ and $G(t,x)=g(t+rx)$, $(t,x)\in \Om_T$.
 Then $G\in C([0,T];H^s(\Om))$.
\end{prop}
 Let us prove Lemma \ref{l.pF1}.
 Observe that, by definition of $\phi(\cdot;\epsl,\tilde\mu)$, and thanks to\refq{v.F1},
 the support of $F_{\epsl,\tilde\mu}|_{\Om_{\tilde\mu}}$ is a subset of the set
$$
 E(\epsl,\tilde\mu)= K_{\epsl}(\tilde t) \cup  (\ove{\Om_{\tilde t+\epsl}}\cap \ove{D^C})
  \cup  (\ove{\Om_{\tilde t+\epsl}}\cap \ove D \cap \{b-\epsl\le x\le b \}).
$$
 Firstly, let $(t,x)\in K_{\epsl}(\tilde t) \cup (\ove{\Om_{\tilde t+\epsl}}\cap \ove{D^C})$.
 Then we have $t-x\le \tilde t+\epsl$, and so
$$ \xi(\mu^{-1}(t-x)) <\xi(\mu^{-1}(\tilde t+\epsl)) < \xi(\mu^{-1}(\tilde \mu-\delta+\epsl)), $$
 since  the functions $\xi$ and $\mu^{-1}$are smooth and non decreasing, and
 $\delta < a(\tilde t) = \tilde\mu-\tilde t$.
 So, for $\epsl$ sufficiently small and some $c>0$ (values that are independent of $t,x$),
 we have
\begin{equation}
\label{m.t-x}
 \xi(\mu^{-1}(t-x)) <\tilde\xi- c\epsl,\quad
 (t,x)\in K_{\epsl}(\tilde t) \cup (\ove{\Om_{\tilde t+\epsl}}\cap \ove{D^C}). 
\end{equation}
 Secondly, let $(t,x)\in \ove{\Om_{\tilde t+\epsl}}\cap \ove D \cap \{b-\epsl\le x\le b \}$.
 Then $t-\frac{x}k \le \nu(t)- \frac{\delta-\epsl}k $ and so, for $\epsl$ sufficiently small
 and some $c>0$,
$$
 \xi(\nu^{-1}(t-\frac{x}k)) \le  \xi(\nu^{-1}(\nu(t)- \frac{\delta-\epsl}k))
 < \tilde\xi -C\epsl.
$$ 
 We thus have
\begin{equation}
\label{m.t-xk}
\xi(\nu^{-1}(t-\frac{x}k))  < \tilde\xi -C\epsl, \quad
  (t,x)\in \ove{\Om_{\tilde t+\epsl}}\cap \ove D \cap \{b-\epsl\le x\le b \}.
\end{equation}
 Since $F_1$ is expressed in terms of $f'_2(t-x)$, $f_2(t-x)$ in $D^C$, and
 in terms of $f'_3(t-\frac{x}k)$, $f_3(t-\frac{x}k)$ in $D$, and since the support
 of $F_{\epsl,\tilde\mu}$ is contained in $E(\epsl,\tilde\mu)$, then, thanks
 to\refq{m.t-x},\refq{m.t-xk}, we see that $F_{\epsl,\tilde\mu}$ can be expressed
 in terms of $f|_{(-\infty,r)}$ and $f'|_{(-\infty,r)}$, $r=\tilde\xi-c\epsl$ only.
 Hence, thanks to Proposition \ref{p.fF}, the conclusion follows.
 \qed
\section{Proof of the main results}
\label{s.3}
\subsection{Proof of Corollary \ref{c.main}} \ \\
\label{ss.21}
 Firstly, notice that $\alpha(t)\neq 0 \egal \dot a(t)\neq \frac{k}{1+k}$.\\
\noindent 1) If $T\le \mu_0$ then $g=0$ in $(0,T)$, and if $T> \mu_0$ then
 $g\neq 0$ since $g|_{(\mu_0,T)} \not\in H^{r_0(1-s^\ast/T)-1}(\mu_0,T)$.
 Hence, the knowledge of $g$ provides $T\le \mu_0$ or $T>\mu_0$.\\
2)
 \begin{itemize}
\item Let $\mu\in [\mu_0,T]$. Thanks to Theorem \ref{t.1}, we can construct
$$ \xi = \inf\{r>0; \; g|_{(0,\mu)} \in H^{r_0(1-r/T)-1}(0,\mu)\}, $$
 and so the invertible function $\mu\mapsto \xi$ from $[\mu_0,T]$ into $[0,s^\ast]$.
 (This implies that $s^\ast$ is recovered too).
 Putting $t=\frac12(\mu+\xi)$, we recover $t_{s^*}$ which is $t$ for $\mu=T$, and also
 the functions $t\mapsto \xi=\xi(t)$, $t\mapsto \mu(t)$,  $t\mapsto a(t)=\frac12(\mu(t)-\xi(t))$,
 for $t\in [t_0,t_{s^*}]$.
 We then construct the functions $t_{(\cdot)}=(\xi(\cdot))^{-1}$, 
 $t^*(\cdot)=2t_{(\cdot)}-\id$.
\item
 Thanks to the above point and to (i) of Theorem \ref{t.1}, the smooth function
 $\alpha(\cdot)$ can be recover as the unique one such that $\mu\mapsto g(\mu)-\alpha(t)f'(\xi)$
 belongs to $H^{\epsl+r_0(1-\xi/T)}(0,\mu)$ for some $\epsl>0$ and all $\mu\in (0,T)$.
 Then, $k$ is root of the following equation:
\begin{equation}
\label{e.k}
 (\alpha+1+\dot a(\alpha-1))k^2 +(\alpha-1)k + \dot a(1-\alpha)= 0.
\end{equation}
 Denote by $k_1,k_2$ the roots, such that $k_1\le k_2$. We show that $k_1\le 0$.
 A short computation shows that
$$
 (\alpha+1+\dot a(\alpha-1)) = \frac2D\left( \frac{(1-\dot a)^2}{1+\dot a}\right)>0,
 \quad D = k(1+\dot a)+ 1-\dot a/k>0.
$$
%
% Let us prove that $k_1\le 0$.
 
 We have
\begin{equation}
\label{e.k1k2}
 k_1k_2 = \frac{\dot a(1-\alpha)}{\alpha+1+\dot a(\alpha-1)} =\dot a (k_1+k_2).
\end{equation}

 If $\dot a\le 0$ then, the second  equality in\refq{e.k1k2} implies that it is impossible
 to have $0<k_1\le k_2$.

% ??????????
%  If $\dot a>0$ then the first equality in\refq{e.k1k2} implies that it is impossible
% to have $k_1k_2>0$, since $1\pm\alpha>0$.\\
% Hence, $k$ can be defined as the unique positive root of\refq{e.k}.
 
?????????, 
 
\end{itemize}
\qed
\begin{remark}
 Theorem \ref{t.1} allows us to recover $t^\ast(\cdot)=\mu\circ \xi^{-1}$ as:
$$
 t^\ast(s) := \sup\{t>s; \:g|_{[s,t]}\in H^{r_0(1-t/T)-1}([0,t]) \},
$$
 and shows that
$$  t^\ast(s) = \sup\{t>s; \: g_A|_{[s,t]}\in H^{r_0(1-t/T)-1}([0,t]) \} .$$
%
% See also Lemma \ref{l.ga} later.
\end{remark}
%
%%%%%%%%%%%%%%%%%%%%%%%%%%%%%%
%
\subsection{Analysis of the error}
\label{ss.23}
 Let $(u_0,u_1,f)$, $r_0$ as in Theorem \ref{t.1}. Put $u=\tilde P(u_0,u_1,f)$,
 $g=\tilde Z(u_0,u_1,f)$, $u_A=U_A(f)$ and
$$
 u_E=u-u_A,\quad F_A=T_Af, \quad g_A=\partial_x u_A|_{x=0},
 \quad g_E=g-g_A=\partial_x u_E|_{x=0}, 
$$
 where $u_A$ is defined in Section \ref{s.Ansatz}.
 Let us prove the estimate\refq{e.gE} (see  (iv) of Theorem \ref{t.1}).
 For the sake of clarity, we replace $\mu$, $t$, $\xi$, respectively by $\tilde\mu$,
 $\tilde t=\mu^{-1}(\tilde\mu)$, $\tilde\xi=\xi(\tilde t)$.
 Put $u_{E,0}=u_0-u_A(0)$, $u_{E,1}=u_1- \partial_t u_A \big|_{t=0}$.
 In view of Subsection \ref{s.Ansatz}, the function $u_E$ satisfies
\begin{eqnarray}
\label{eq.uEA}
 \left\{ \begin{array}{rll}
  \calL_\gamma u_E &=& -F_A \quad {\rm in} \quad \Om_{\tilde\mu},\\
  u_E|_{x=0,b} &=& 0  \quad {\rm on} \quad (0,\tilde\mu), \\
  u_E\big|_{t=0} &=& u_{E,1}  \quad {\rm on} \quad \Om,\\
  \partial_t u_E\big|_{t=0} &=& u_{E,1} \quad {\rm on} \quad \Om.
\end{array} \right.
\end{eqnarray}
 So we have $u_E=P(u_{E,0},u_{E,1},-F_A)$.
 Recall that, thanks to Lemma \ref{l.T0}, we have $T_0(f)\in L^2(0,\tilde\mu;H^{-s}(\Om))$,
 for all $s>\frac12$. Thanks to Proposition \ref{p.Zinter}, we have
\begin{equation}
\label{ev.ZuE1}
 Z(0,0,T_0(f))\big|_{(0,\tilde\mu)}\in H^{-s}(0,\tilde\mu),\quad \forall s>\frac12.
\end{equation} 
 Let us prove that $u_{E,0}\in  H^{r_0}(\Om)$, $u_{E,1}\in  H^{r_0-1}(\Om)$.
 Observe that $u_A(0)(x)= (f(-x)+f_2(x)+f_2(-x)\Phi_\epsl(x))\chi_{x<a(0)}+
 f_3(-x/k)\Phi_\epsl(x-b+2\epsl))\chi_{x>a(0)}$.
 For $x<a(0)=t_0$ we have
$$ \xi(\mu^{-1}(x))<\xi(\mu^{-1}(t_0))<\xi(\mu^{-1}(\mu_0))= \xi(t_0)=0, $$
 and, similarly, $\xi(\mu^{-1}(-x))\le \xi(\mu^{-1}(0))<0$.
 For $x>a(0)$ we  have
$$ \xi(\nu^{-1}(-x/k))<\xi(\nu^{-1}(-t_0/k))<\xi(\nu^{-1}(\nu_0))=0 .$$
 Hence, $u_A(0)$ can be expressed in terms of $f(\xi)$ for $\xi<0$.
 Since $f|_{(-\infty,0]}\in H^{r_0}(-\infty,0)$, then $u_A(0)\in  H^{r_0}(\Om)$.
 Thanks to the asumption on $u_0$, we then have
 $u_{E,0}\in  H^{r_0}(\Om)$. Similarly, we have  $u_{E,1}\in  H^{r_0-1}(\Om)$.
 Thanks to\refq{v.f2}, the regularity of $f_2|_{(0,\tilde\mu)}$ is given by
 those of $f|_{(0,\tilde\xi)}$, that is, $f_2|_{(0,\tilde\mu)}\in H^{r_0(1-\xi'/T)}((0,\tilde\mu))$,
 for all $\xi'>\tilde\xi$. Thus, thanks to Proposition \ref{p.Zinter}, we have
\begin{equation}
\label{ev.ZuE}
 Z(u_{E,0},u_{E,1},0)\big|_{(0,\tilde\mu)}\in H^{r_0-1}(0,\tilde\mu).
\end{equation} 
 Thanks to Lemma  \ref{l.UC} with $t_1$ replaced by $\tilde t$ and $T$ by $\tilde \mu$,
 and to Lemma \ref{l.suppFF1}, we have
\begin{equation}
\label{r.ZF1F}
 Z(0,0,-F_1)\big|_{(0,\tilde\mu)}=Z(0,0,-F_{\epsl,\tilde\mu})\big|_{(0,\tilde\mu)} .
\end{equation} 
 Thanks to Lemma \ref{l.pF1}, if $\epsl>0$ is sufficiently small, we have
$$
 F_{\epsl,\tilde\mu}\big|_{\Om_{\tilde\mu}} \in  L^2([0,\tilde\mu];
 H^{r_0(1-\tilde\xi/T)+c\epsl-1}(\Om)),
$$
 and so, thanks to\refq{r.ZF1F} and by applying Proposition \ref{p.Zinter}, we obtain
\begin{equation}
\label{ev.ZF1}
 Z(0,0,-F_1)\big|_{(0,\tilde\mu)}\in H^{r_0(1-\tilde\xi/T)+\epsl-1}(0,\tilde\mu),
\end{equation} 
 for some $\epsl>0$ (independent of $\tilde\mu$).
 
 Thanks to\refq{ev.ZuE1},\refq{ev.ZuE}\refq{ev.ZF1}, and since
 $g_E=  Z(u_{E,0},u_{E,1},0)+  Z(0,0,T(0)f)+ Z(0,0,-F_1)$, the proof of\refq{e.gE} is done.

\qed
%%%%%%%%%%%%%%%%%%%%%%%%%%%%%%%%%%%%%%%%%%
\newpage

\section{Appendix: the function G}
 Let $I=(0,1)$ and a dense sequence $\{a_n\}_{n\in\N^*}$ in $\ove I$.
 We set
$$ f_n(x)=((x-a_n)_+)^{1/2-a_n}, $$
$$ G(x) = \sum_{n\in\N^*} \frac1{2^n} f_n(x), \quad x\in I, $$
 where $z_+=\max(0,z)$ for $z\in\R$. The function $G$ is increasing.

\vskip1em

  For $0<s<1$ we set the following Sobolev space:
$$ H^s(I) = \left\{q\in L^2(I);\; \int\int_{I\times I} \frac{|q(x)-q(y)|^2}{|x-y|^{1+2s}} \; dx\, dy <\infty\right\} .$$
\begin{lemma}
\label{l.faHs}
 Let $b\in (0,1]$, $r> -\frac12$, $s\in (0,1)$, $a\in [0,b)$.
 Set $f(x) = ((x-a)_+)^r$, $I_b=(0,b)$.
 We have $f \in H^s(I_b)$ if, and only if, $r>s-1/2$. In such a case, we have
\begin{equation}
\label{m.ff}
 \int\int_{I\times I} \frac{|f(x)-f(y)|^2}{|x-y|^{1+2s}} \; dx\, dy \le
 C_s \left(\frac1{2r+1} + \frac{r^2}{2r-2s+1}\right) (b-a)^{2r-2s+1},
\end{equation}
 for some $C_s>0$.
\end{lemma}
 Proof. Firstly, let $b=1$. We have
\begin{eqnarray*}
 J&:=& \int\int_{I_1\times I_1} \frac{|f(x)-f(y)|^2}{|x-y|^{1+2s}} \; dx\, dy =
  2\int_0^1 dy\left(  \int_0^y\frac{|f(x)-f(y)|^2}{|x-y|^{1+2s}} \; dx\right) \\
 &=& 2(K_1+K_2),\\
 K_1 &:=& \int_a^1 dy\left(\int_0^a \frac{(y-a)^{2r}}{(y-x)^{1+2s}} \; dx\right) ,\\
 K_2 &:=& \int_a^1 dy\left(\int_a^y \frac{((y-a)^r-(x-a)^r)^2}{(y-x)^{1+2s}} \; dx\right).
\end{eqnarray*}
 We have
\begin{eqnarray*}
 K_1 &=&  \frac1{2s}\int_a^1  (y-a)^{2r}\left[\frac1{(y-x)^{2s}} \right]_0^a dy 
 =  \frac1{2s}\int_a^1 \left( (y-a)^{2r-2s}-\frac{ (y-a)^{2r}}{y^{2s}} \right)dy.
\end{eqnarray*}
 If $a=0$, then $K_1=0$. If $a>0$, then $K_1<\infty$ if, and only if, $2r>2s-1$.
 In such a case, we have
\begin{equation}
\label{m.K1}
 K_1 \le  \frac1{2s(2r-2s+1)} (1-a)^{2r-2s+1}.
\end{equation}
 Let $2r>2s-1$. We have
\begin{eqnarray*}
 K_2 &=& \int_0^{1-a} dy\left(\int_0^y \frac{(y^r-x^r)^2}{(y-x)^{1+2s}} \; dx\right)
 = \int_0^{1-a} y^{2r-2s}dy\left(\int_0^1 \frac{(1-t^r)^2}{(1-t)^{1+2s}} \; dt\right)\\
 &=& \frac{C(r,s)}{2r-2s+1}(1-a)^{2r-2s+1},
\end{eqnarray*}
 where
\begin{eqnarray}
\nonumber
 C(r,s) &=& \int_0^1 \frac{(1-t^r)^2}{(1-t)^{1+2s}} \; dt =
 \int_0^{1/2} \frac{(1-t^r)^2}{(1-t)^{1+2s}} \; dt + \int_{1/2}^1 \frac{(1-t^r)^2}{(1-t)^{1+2s}} \; dt \\
\label{m.K2}
 &\le & C_s(\frac1{2r+1} + \frac{r^2}{2r-2s+1}).
\end{eqnarray}
 Since $C(r,s)>0$, then $K_2=+\infty$ if $2r\le2s-1$.
 Hence, the sum $K_1+K_2$ converges iff $2r>2s-1$.
 If $2r>2s-1$, thanks to\refq{m.K1} and\refq{m.K2}, we obtain\refq{m.ff}.

 Secondly, the case $b\in (0,1)$ is easily proved by setting $a=a'b$, $x=x'b$, $y=y'b$.
\qed
\begin{lemma}
 For $0<s< 1$ and $b\in (0,1]$, we have $G\in H^s(0,b)$ if $s<1-b$ and $G\not\in H^s(0,b)$ if $s>1-b$.
\end{lemma}
 Proof. For $x,y\in I$, we have, thanks to the Schwarz inequality,
\begin{equation}
\label{i.Schw}
 |G(x)-G(y)|^2 \le \left(\sum_{n\ge 1} \frac1{2^n}\right) \left(\sum_{n\ge 1} \frac1{2^n} |f_n(x)-f_n(y)|^2\right)
 =  \sum_{n\ge 1} \frac1{2^n} |f_n(x)-f_n(y)|^2.
\end{equation}
 Let $I_b=(0,b)$, $A_b=\{n\in \N^*; \; a_n\ge b\}$, %note: n\in A_b si a_n=1
 $B_b=\N^*\setminus A_r = \{n; a_n <b \}$.\\ %Note:  n\in B_b=> a_n\neq 1
 For all $n\in B_b$, thanks to Lemma \ref{l.faHs}, we have $f_n\in H^{1-b}(0,1)$, since $1/2-a_n>(1-b)-1/2$.
 For all $n\in A_b$, we have $f_n\in H^{1-b}(I_b)$, since $f_n|_{I_b}=0$.\\
 Let $0<s<1-b$. By using\refq{i.Schw}, and\refq{m.ff}, we have
\begin{eqnarray*}
 J_{b,s} &:=& \int\int_{I_b\times I_b} \frac{|G(x)-G(y)|^2}{|x-y|^{1+2s}} \; dx\, dy \le
  \sum_{n\in B_b} \frac1{2^n}\int\int_{I_b\times I_b} \frac{|f_n(x)-f_n(y)|^2}{|x-y|^{1+2s}} \; dx\, dy \\
 & \le & C_s \sum_{n\in B_b} \frac1{2^n} (\frac1{1-a_n}+ \frac1{1-a_n-s}) (b-a_n)^{2(1-s-a_n)} \\
 &\le & C_s \sum_{n\in B_b} \frac1{2^n} (\frac1{1-b}+ \frac1{1-b-s}) (b-a_n)^{2(1-s-a_n)}
 <\infty
\end{eqnarray*}
 since $(b-a_n)^{2(1-s-a_n)}\le 1$ for all $n\in B_b$, $0<s<1-b$.

 Let $s\in (1-b,1)$. For all $n\in\N^*$ and $x>y$ we have $G(x)-G(y) \ge f_n(x)-f_n(y)$.
 Fix $n\in A_{1-s}\cap B_b$, that is, $1-s\le a_n<b$.
 Thanks to Lemma \ref{l.faHs}, we have $f_n\not\in H^s(I_b)$, and then
\begin{eqnarray*}
 J_{b,s} &\ge & \frac1{2^n}\int\int_{I_b\times I_b} \frac{|f_n(x)-f_n(y)|^2}{|x-y|^{1+2s}} \; dx\, dy 
 = \infty.
\end{eqnarray*}
 This ends the proof.
\qed

\section{Proof of Theorem \ref{t.L1}}
 Let $F\in L^2(\Om_T)$, $v_0\in H^1_0(\Om)$, $v_1\in L^2(\Om)$.
 Denote $M^1:=\{v\in C([0,T];H^1_0(\Om))$, $\partial_t v\in C([0,T];L^2(\Om))\}$,
 $M^1_0=\{v\in M$; $v|_{t=0}=0$, $\partial_t v|_{t=0}=0\}$,
\subsection{ Energy estimate.}
 Put
$$ E(t)(v)= \frac12\int_\Om |\partial_t v|^2 + \frac12\int_\Om \gamma(t,\cdot)|\partial_x v|^2, \quad v\in M^1.$$
 We claim that, for all $v\in M^1$ such that $L_\gamma v=:f \in L^2(\Om_T)+W$,
  the following (standart) estimate, which implies\refq{in.en1}, holds.
\begin{equation}
\label{in.en}
 E(t)(v) \le C \left(\|f\|^2_{L^2(0,t;\Om)}+E(0)(v) \right),\quad \forall t\in [0,T],
\end{equation}
 for some constant $C$. \\ %=C(\delta,\rho)$.\\
 Proof. It is sufficient to show\refq{in.en} for $t=T$. Assume that $f\in L^2(\Om_T)$.
 Put $\rho=\sup_Q \frac{|\dot\gamma|}{\gamma}$ and $\Pi_0\in C^1([0,T];(0,+\infty))$
 such that $\delta^{-1}\Pi_0\le -\Pi'_0$ for some  $\delta\in (0,\frac1{\rho})$.
 (For example, $\Pi_0=e^{-\frac{t}{\delta}}$).
 Put
$$ Q(v) = \int_0^T E(t)(v) \; \Pi_0 dt,\quad  C_0(f) = \int_Q f^2 \Pi_0 .$$
 We formally have, thanks te the Schwarz inequality,
\begin{eqnarray*}
 \delta^{-1}Q(v) &\le & -\int_0^T E(t)(v)  \; \Pi'_0 dt =  [-E(t)(v) \; \Pi_0(t)]_0^T -\int_0^T \frac{dE(t)(v)}{dt}  \; \Pi_0 dt\\
 &\le & E(0)(v)\,\Pi_0(0)-E(T)(v) \; \Pi_0(T) - \frac12\int_Q \Pi_0 \dot\gamma |\partial_t v|^2 - \int_Q\Pi_0 f \partial_t v\\
 &\le & E(0)(v)\,\Pi_0(0)-E(T)(v) \; \Pi_0(T)+  \rho Q(v) + \sqrt{2C_0(f)} \sqrt{Q(v)},
\end{eqnarray*}
 Hence, we obtain
$$ (\delta^{-1}-\rho)Q(v) +E(T)(v) \; \Pi_0(T)  \le  E(0)(v)\,\Pi_0(0)+ \sqrt{2C_0(f)}\sqrt{Q(v)},$$
 and so,
\begin{equation}
\label{in.en2}
 Q(v)+E(T)(v) \; \Pi_0(T) \le C (C_0(f)+E(0)(v)\,\Pi_0(0)).
\end{equation}
 Then\refq{in.en2} follows.
\subsection{Uniqueness}
 Consequently, if $v\in M^1_0$ satisfies\refq{eq.v} with $F=0$, then $E(t)(v)\equiv 0$ for all $t$,
 and so $v\equiv 0$. This shows that Problem\refq{eq.v} admits at most one solution in $M^1$.
\subsection{Existence}
 Let $(\lambda_j,e_j)_{1\le j}$ be the familly of spectral values of the positive
 operator $-\Delta_x$ in $H^1_0(\Om)$, i.e such that $(e_i,e_j)_{L^2(\Om)}=\delta_{ij}$,
 $-\Delta e_j=\lambda_j e_j$, and $\lambda_j \nearrow +\infty$.
 The data $v_0,v_1$, $F$ are then written
 $v_0=\sum_{j=1}^\infty v_{0,j} e_j$, $v_1=\sum_{j=1}^\infty v_{1,j} e_j$, $F(t,\cdot)=\sum_{j=1}^\infty F_j(t) e_j$, 
 with
 $$ \sum_{j=1}^\infty\{\lambda_j |v_{0,j}|^2 + |v_{1,j}|^2 + \int_0^T |F_j(t)|^2 dt\} <\infty.$$
 Let $N\in \N^*$, and put $E_N={\rm span}\{e_1,\ldots,e_N\}$,
 $V_{k,N}=(v_{k,1},\ldots,v_{k,N})$, $k=0,1$, $F_N=\sum_{j=1}^N F_j(t) e_j$,
 $B_N(t)=(b_{i,j}(t))_{1\le i,j\le N}$ with $b_{i,j}(t)=(\nabla e_i,\nabla e_j)_{L^2(\Om;\gamma(t,\cdot)dx)}$,
 and consider the following vectorial differential equation:
 find $V_N(t)=(v_1(t),\ldots,v_N(t))$ such that
$$ \frac{d^2}{dt^2} V_N(t)  +  V_N(t)B_N(t) = F_N(t),\quad 0\le t\le T,$$
 with the initial condition $V_N(0)= V_{0,N}$, $\frac{d}{dt} V_N(0)= V_{1,N}$.
 Since $B_N(\cdot)$ is continuous, the theorem of Cauchy-Lipschitz implies existence
 and uniqueness for $V_N(t)$.
 Note that $B_N(t)$ is positive since, for all $U=(u_1,\ldots,u_N)$, setting
 $u(x)=\sum_{j=1}^N u_j e_j(x)$, we have
$$
 UB_N(t)\,^t\!U = \int_\Om |\nabla_x u|^2\gamma(t,x) dx \ge C\|\nabla_x u\|^2_{L^2(\Om)}
 = C  \sum_{j=1}^N \lambda_j |u_j|^2,
 $$
 where $C$ is a constant such that $0<C\le \gamma$ in $Q$.
 Let $v_N(t)=\sum_{j=1}^N v_j(t) e_j(x)$.
 Then, a standart energy estimate for $E_N(t)(v_N) = \frac12 (\dot V_N^2(t) + V_N(t) B_N(t)V_N(t))$,
 as above, implies that there exists a positive constant $C$ such that
$$
 \|\dot v_N(t)\|_{L^2(\Om)} + \|\partial_x v_N(t)\|_{L^2(\Om)} \le 
 C (\|v_0\|_{H^1(\Om)} +\|v_1\|_{L^2(\Om)}+ \|F\|_{L^2(\Om)}),\quad 0\le t\le T.
$$
 Passing to the limit $N\to +\infty$, we can conclude by standard arguments
 that $(v_N)_N$ converges to a function $v\in C([0,T];H^1_0(\Om))$ satisfying\refq{eq.v}.

 The proof of Theorem \ref{t.L1} in done in the case $F\in L^2(\Om_T)$.
 The case $F\in W$ is similar.

\bibliographystyle{../../MACROS/iopart_BibTeX/iopart-num}

\begin{thebibliography}{99}
\bibitem{DAU.ANA} R. Dautray. and J.L. Lions, {\it Analyse math\'ematique et calcul num\'erique
 pour les sciences et les techniques},{\bf 7}, Masson, Paris Milan Barcelone Mexico (1988).

\bibitem{GRI.SIN} P. Grisvard, {\it Singularities in boundary value problems}, RMA{\bf 22}, Masson, Springer-Verlag (1992).
\end{thebibliography}

\end{document}